\documentclass [11pt] {amsart}
\usepackage{bbm}
\usepackage{amssymb}
\usepackage{epic}
\usepackage{ecltree}
\usepackage{tikz}
\newtheorem{theorem}{Theorem}
\newtheorem{lemma}[theorem]{Lemma}

\newtheorem*{theoremA}{Theorem A}
 \newtheorem*{theoremB}{Theorem B}
\theoremstyle{definition}
\newtheorem{definition}[theorem]{Definition}

\newtheorem{remark}[theorem]{Remark}

\begin{document}

\title{Metric embeddings of Laakso graphs into Banach spaces}
 \author{S.\ J.\ Dilworth, Denka  Kutzarova and Svetozar Stankov}

\address{Department of Mathematics, University of South Carolina, Columbia, SC 29208, U.S.A.}
\email{dilworth@math.sc.edu}

\address{Institute of Mathematics, Bulgarian Academy of
Sciences, Bulgaria} \curraddr{Department of Mathematics, University of Illinois
at Urbana-Champaign, Urbana, IL 61801, U.S.A.} \email{denka@math.uiuc.edu}

\email{erejnion@gmail.com}

\begin{abstract} Let $X$ be  Banach  space which is not super-reflexive. Then, for each $n\ge1$ and $\varepsilon>0$, we exhibit  metric embeddings of  the Laakso graph $\mathcal{L}_n$   into $X$ with distortion less than $2+\varepsilon$ and  into $L_1[0,1]$ with distortion $4/3$. The distortion of an embedding of $\mathcal{L}_2$ (respectively, the diamond graph $D_2$) into $L_1[0,1]$ is at least $9/8$ (respectively, $5/4$).
\end{abstract} 
\thanks{S. J. Dilworth was supported by Simons Foundation Collaboration Grant No. 849142.
 Denka Kutzarova was supported by Simons Foundation Collaboration Grant No. 636954.}

\maketitle
\section{Introduction}

 James \cite{J1972} introduced the important property of super-reflexivity: a Banach space $X$ is  super-reflexive  if every Banach space $Y$ which is finitely representable in $X$ is reflexive.  Enflo  \cite{E} showed that super-reflexivity of $X$  is equivalent to $X$ having an equivalent uniformly convex norm.  

Let us recall the definition of the diamond and Laakso graphs.

\begin{definition}\label{D:Diamonds}
The diamond graph of level $0$ has two vertices joined by an
edge of length $1$ and is denoted by $D_0$. The {\it  diamond
graph} $D_n$ is obtained from $D_{n-1}$ in the following way.
Each edge $uv$ of  $D_{n-1}$ is replaced by a
quadrilateral $u, a, v, b$, with edges $ua$, $av$, $vb$, $bu$ of length $1$.
(See Figure \ref{F:Diamond2}.)
\end{definition}

Definition \ref{D:Diamonds} was  introduced in
\cite{GNRS04}.

\begin{figure}
\begin{center} {
\begin{tikzpicture}[scale=0.2]

\node (A)  at (16,0) {$\bullet$};
 \node (B)  at (16,8) {$\bullet$};
 \node (C)  at (8,16) {$\bullet$};
 \node (D)  at (24,16) {$\bullet$};
 \node (E)  at (16,24) {$\bullet$};
 \node (F)  at (16,32) {$\bullet$};

\node (A1)  at (-8,0) {$\bullet$};
 \node (B1)  at (-8,8) {$\bullet$};
 \node (C1)  at (-16,16) {$\bullet$};
 \node (D1)  at (0,16) {$\bullet$};
 \node (E1)  at (-8,24) {$\bullet$};
 \node (F1)  at (-8,32) {$\bullet$};

\draw (A1) edge  (B1);
\draw (B1) edge    (C1);
\draw (B1) edge   (D1);
\draw (C1) edge   (E1);
\draw (D1) edge  (E1);
\draw (E1) edge    (F1);

\node (b1)  at (14,10) {$\bullet$};
 \node (b2)  at (11,11) {$\bullet$};
 \node (b3)  at (13,13) {$\bullet$};
 \node (b4)  at (10,14) {$\bullet$};
\draw (b1) edge   (B);
\draw (b1) edge    (b2);
\draw (b1) edge    (b3);
\draw (b2) edge   (b4);
\draw (b3) edge  (b4);
\draw (C) edge  (b4);

\node (g1)  at (22,18) {$\bullet$};
 \node (g2)  at (19,19) {$\bullet$};
 \node (g3)  at (21,21) {$\bullet$};
 \node (g4)  at (18,22) {$\bullet$};
\draw (g1) edge  (D);
\draw (g1) edge  (g2);
\draw (g1) edge   (g3);
\draw (g2) edge   (g4);
\draw (g3) edge   (g4);
\draw (E) edge  (g4);

 \node (a1)  at (16,2) {$\bullet$};
 \node (a2)  at (15,4) {$\bullet$};
 \node (a3)  at (17,4) {$\bullet$};
 \node (a4)  at (16,6) {$\bullet$};
\draw (a1) edge   (A);
\draw (a2) edge  (a1);
\draw (a3) edge  (a1);
\draw (a4) edge  (a2);
\draw (a4) edge  (a3);
\draw (B) edge  (a4);

 \node (e1)  at (16,26) {$\bullet$};
 \node (e2)  at (15,28) {$\bullet$};
 \node (e3)  at (17,28) {$\bullet$};
 \node (e4)  at (16,30) {$\bullet$};
\draw (e1) edge  (E);
\draw (e2) edge   (e1);
\draw (e3) edge  (e1);
\draw (e4) edge (e2);
\draw (e4) edge   (e3);
\draw (F) edge   (e4);

 \node (d1)  at (18,10) {$\bullet$};
 \node (d2)  at (21,11) {$\bullet$};
 \node (d3)  at (19,13) {$\bullet$};
 \node (d4)  at (22,14) {$\bullet$};
\draw (d1) edge  (B);
\draw (d1) edge   (d2);
\draw (d1) edge  (d3);
\draw (d2) edge   (d4);
\draw (d3) edge   (d4);
\draw (D) edge   (d4);

 \node (f1)  at (10,18) {$\bullet$};
 \node (f2)  at (13,19) {$\bullet$};
 \node (f3)  at (11,21) {$\bullet$};
 \node (f4)  at (14,22) {$\bullet$};
\draw (f1) edge   (C);
\draw (f1) edge   (f2);
\draw (f1) edge   (f3);
\draw (f2) edge   (f4);
\draw (f3) edge  (f4);
\draw (E) edge   (f4);

\end{tikzpicture}} \caption{The Laakso graphs $\mathcal{L}_1$ and $\mathcal{L}_2$} \label{fig: L1L2}
\end{center}
 \end{figure}
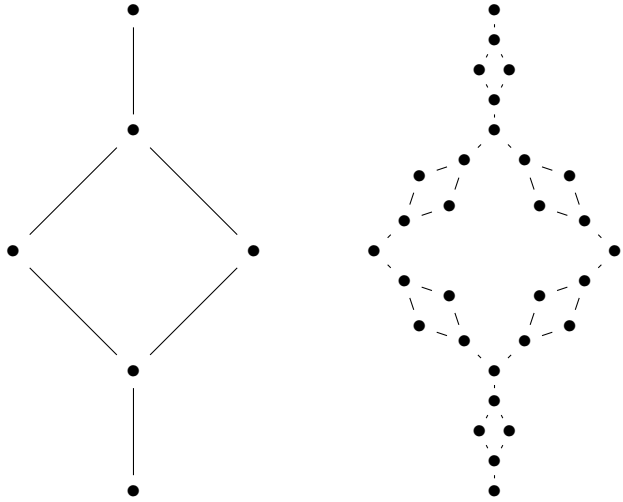

\begin{definition}\label{D:Laakso} The 
Laakso graph of level $0$ has two vertices joined by an edge of
length $1$ and is denoted $\mathcal{L}_0$. The {\it Laakso graph}
$\mathcal{L}_n$ is obtained from $\mathcal{L}_{n-1}$ according to
the following procedure. Each edge $uv\in E(\mathcal{L}_{n-1})$ is
replaced by the graph $\mathcal{L}_1$ exhibited in Figure
\ref{fig: L1L2} in which each edge has length $1$.
\end{definition}

Definition \ref{D:Laakso} was introduced in \cite{LP} based on an
idea of Laakso \cite{Laa00}. 

Let $f \colon (M,\rho) \rightarrow (N,\sigma)$ be a bilipschitz mapping beween metric spaces. The distortion of $f$ is defined to be the infimum of $b/a$,
where $a,b$ are positive constants such that $$a \rho(x,y) \le \sigma(f(x),f(y)) \le b \rho(x,y) \qquad (x,y \in M).$$

Bourgain \cite{B1986} characterized Banach spaces which are not  super-reflexive  as those for which the binary trees $B_n$   of depth $n$ embed with uniformly bounded distortion. Subsequently, Johnson and Schechtman \cite{JS} characterized Banach  spaces which are not  super-reflexive   as those for which the diamond graphs $D_n$ and the Laakso graphs $\mathcal{L}_n$ embed with uniformly bounded distortion. The best known   estimate in the literature for the distortion of embeddings of $D_n$  into Banach spaces which are not super-reflexive, due to Pisier \cite{P}, 
is  $2+\varepsilon$ for every $\varepsilon>0$, while the best known estimate for the distortion of embeddings of $D_n$ into $L_1[0,1]$, due to Lee and Rhagavendra \cite {LR}, is $4/3$.

In the present article we construct embeddings of $\mathcal{L}_n$ into arbitrary  Banach spaces which are not super-reflexive with disortion $2 + \varepsilon$ and into $L_1[0,1]$ with distortion $4/3$. We also show that $\mathcal{L}_2$ does not embed into $L_1[0,1]$ with distortion smaller than $9/8$.

\begin{figure}
\begin{center}
{
\begin{tikzpicture}
  [scale=.15,auto=left,every node/.style={circle,draw}]
  \node (n1) at (16,0) {\hbox{~~~}};
  \node (n2) at (5,5)  {\hbox{~~~}};
  \node (n3) at (11,11)  {\hbox{~~~}};
  \node (n4) at (0,16) {\hbox{~~~}};
  \node (n5) at (5,27)  {\hbox{~~~}};
  \node (n6) at (11,21)  {\hbox{~~~}};
  \node (n7) at (16,32) {\hbox{~~~}};
  \node (n8) at (21,21)  {\hbox{~~~}};
  \node (n9) at (27,27)  {\hbox{~~~}};
  \node (n10) at (32,16) {\hbox{~~~}};
  \node (n11) at (21,11)  {\hbox{~~~}};
  \node (n12) at (27,5)  {\hbox{~~~}};

  \foreach \from/\to in {n1/n2,n1/n3,n2/n4,n3/n4,n4/n5,n4/n6,n6/n7,n5/n7,n7/n8,n7/n9,n8/n10,n9/n10,n10/n11,n10/n12,n11/n1,n12/n1}
    \draw (\from) -- (\to);

\end{tikzpicture}
} \caption{The diamond graph $D_2$.}\label{F:Diamond2}
\end{center}
\end{figure}
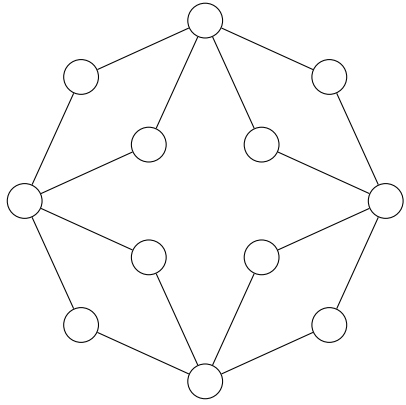

\section{Results} The embeddings of $\mathcal{L}_n$ which we define depend on the following characterization of not being super-reflexive. Its  negation is the characterization of super-reflexivity known as $J$-convexity.
\begin{theoremA} \cite{JS1972, SS} $X$ is not super-reflexive if and only if, for each $m \ge 1$ and $\varepsilon>0$, there exist $e_1,\dots, e_m$ in the unit ball of $X$ such that, for each $1 \le j \le m$, we have \begin{equation} \label{eq: J-convex}
\| e_1+\dots + e_j - e_{j+1}-\dots-e_m\| \ge m - \varepsilon. \end{equation}
\begin{remark} It follows easily from Theorem~A that if $X$ is not super-reflexive then, for each $n\ge1$ and  $\varepsilon>0$,  $B_n$ embeds into $X$ with distortion $1+\varepsilon$. This is not true, however, for $D_n$ and $\mathcal{L}_n$ if $n \ge 2$.
\end{remark} 
\end{theoremA} We wil make use of the following two consequences of Theorem~A.
\begin{lemma} \label{lem: J-convex} Suppose $X$ is not super-reflexive.  Let $(e_i)_{i=1}^m $ be as in Theorem~A.  If  $\max A < \min B$ then
$$\| \sum_{i \in A} e_i - \sum_{i \in B} e_i\| \ge |A| + |B| - \varepsilon.$$
\end{lemma} 
\begin{proof}  This follows  at once from  \eqref{eq: J-convex}  and the triangle inequality. \end{proof}
\begin{lemma} \label{lem: J-convex2}  Suppose $X$ is not super-reflexive.  Let $(e_i)_{i=1}^m $ be as in Theorem~A. If $\max A < \min B$ or $\max B < \min A$ then
$$ \|\sum_{i \in A}  \varepsilon_i e_i + \sum_{ i \in B} e_i\| \ge |B| - \varepsilon.$$ 
for all choices of signs $\varepsilon_i = \pm1$.\end{lemma}
\begin{proof} Let $A^+ = \{ i \in A \colon \varepsilon_i = 1\}$ and let  $A^-= \{ i \in A \colon \varepsilon_i = -1\}$.
If $|A^+| \ge |A^-|$ then \begin{align*}
\|\sum_{i \in A} e_i + \sum_{i \in B} e_i\| &\ge \|\sum_{i \in A^+} e_i + \sum_{i \in B} e_i\| - |A^-|\\
&\ge |A^+| + |B| - \varepsilon - |A^-|\\ \intertext{(by Lemma~\ref{lem: J-convex})}
&\ge |B| - \varepsilon.
\end{align*} On the other hand, if $|A^-| > |A^+|$ then 
\begin{align*} \|\sum_{i \in A} e_i + \sum_{i \in B} e_i\| &\ge \|-\sum_{i \in A^-} e_i + \sum_{i \in B} e_i\| - |A^+|\\
&\ge |A^-| + |B| - \varepsilon - |A^+|\\ \intertext{(by Lemma~\ref{lem: J-convex})} &\ge 1+ |B| - \varepsilon.
\end{align*}
\end{proof} \begin{theorem} \label{thm: superreflexive}
 Suppose $X$ is not super-reflexive. Then, for each $\varepsilon>0$ and  $n \ge 1$, there exists a mapping $f_n \colon \mathcal{L}_n \rightarrow X$ such that, for all $a,b \in \mathcal{L}_n$,  
\begin{equation}  \label{eq: bilipschitz} \frac{1}{2} d(a,b) -\varepsilon \le \|f_n(a)-f_n(b)\| \le d(a,b).
\end{equation} \end{theorem}
\begin{proof} Let $\varepsilon>0$ be fixed.  For each $n \ge 1$, select vectors $(e^n_i)_{i=1}^{4^n}$ satisfying Lemma~\ref{lem: J-convex}
for $m = 4^n$. We define the mappings $f_n$ inductively. 

 We begin with the base case $n=1$. Label the vertices of $\mathcal{L}_1$ as shown in Figure~\ref{figure: L1}.
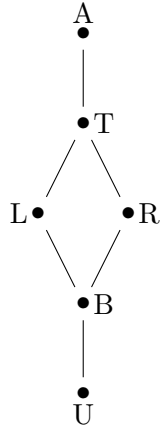
\begin{figure} \label{fig: L} \begin{center} {
\begin{tikzpicture}[scale=0.15]
 \node (A)  at (66,0) {$\bullet$};
\node[below] at (66, 0) {U};
\node (B) at (66,8) {$\bullet$}; \node[right] at (66,8) {B};
\node (C) at (62,16) {$\bullet$}; \node[left] at (62,16) {L};
\node (D) at (70,16) {$\bullet$}; \node[right] at (70,16) {R};
\node (E) at (66,24) {$\bullet$}; \node[right] at (66,24)  {T};
\node (F) at (66,32) {$\bullet$};
\node[above] at (66,32) {A};
\draw (A)  edge (B);
\draw (B)  edge (C); \draw  (B)  edge (D);\draw  (C)  edge (E); \draw  (D)  edge (E);
\draw  (E)  edge (F);
\end{tikzpicture}
}\end{center} \caption{The Laakso graph $\mathcal{L}_1$} \label{figure: L1} \end{figure}
We define $f_1 \colon \mathcal{L}_1 \rightarrow X$  as follows:
$$f_1(A)=0, f_1(T)= e^1_1, f_1(L)= e^1_1+e^1_2, f_1(R)=e^1_1+e^1_3$$$$ f_1(B)=e^1_1+e^1_2+e^1_3, f_1(U)= e^1_1+e^1_2+e^1_3+e^1_4.$$
Using Lemma~\ref{lem: J-convex} it is easily checked that  $f_1$ satisfies, for all $a,b \in \mathcal{L}_1$,
$$d(a,b) - \varepsilon \le  \|f_1(a)-f_1(b)\| \le d(a,b).$$
For example,  $f_1(L) - f_1(R) =  e^1_2 - e^1_3$, so $2-\varepsilon \le \|f_1(L)-f_1(R)\| \le 2$ as required.

Now suppose $n\ge2$. We regard $\mathcal{L}_n$ as being obtained from $\mathcal{L}_1$ by replacing each edge of $\mathcal{L}_1$ by a copy of $\mathcal{L}_{n-1}$. Thus $\mathcal{L}_n$ is  composed of $6$  copies of $\mathcal{L}_{n-1}$, labelled as $Y,C,D,E,F$ and $Z$  in Figure~\ref{figure: Ln}.

We have labelled the vertices $A,T,L,R,B$ and $U$ of $\mathcal{L}_n$ which correspond to the vertices of $\mathcal{L}_1$.  The correspondence between  $\mathcal{L}_{n-1}$ and each of its copies in $\mathcal{L}_n$, namely  $Y,C,D,E,F,$ and $Z$,   is the natural `downward' correspondence in which the vertex $A$ of $\mathcal{L}_{n-1}$ is mapped to the vertices $A,T,T,L,R,$ and $B$ of $\mathcal{L}_n$  respectively.
 Note that   the vertex T  of $\mathcal{L}_n$ corresponds to 
the vertex $U$ of $Y$ and to the vertex $A$ of $C$ and $D$.  There are similar correspondences for $L,R$ and $B$.
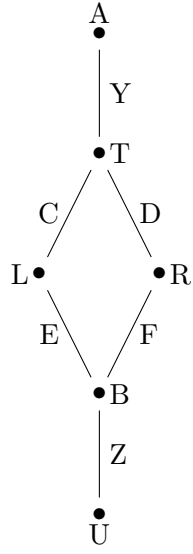
\begin{figure} \label{fig: Ln}\begin{center} {
\begin{tikzpicture}[scale=0.20]
 \node (A)  at (16,0) {$\bullet$};
\node[below] at (16, 0) {U};
\node (B) at (16,8) {$\bullet$}; \node[right] at (16,8) {B};
\node (C) at (12,16) {$\bullet$}; \node[left] at (12,16) {L};
\node (D) at (20,16) {$\bullet$}; \node[right] at (20,16) {R};
\node (E) at (16,24) {$\bullet$}; \node[right] at (16,24)  {T};
\node (F) at (16,32) {$\bullet$};
\node[above] at (16,32) {A};
\draw (A)  edge (B);
\draw (B)  edge (C); \draw  (B)  edge (D);\draw  (C)  edge (E); \draw  (D)  edge (E);
\draw  (E)  edge (F);
   \node [left] at (14,12) {E};  \node [right] at (18,12) {F};  \node [right] at (18,20) {D};
 \node [left] at (14,20) {C};  \node [right] at (16,4) {Z};  \node [right] at (16,28) {Y};
\end{tikzpicture}} \end{center}  \caption{The Laakso graph $\mathcal{L}_n$}\label{figure: Ln} \end{figure}

Let $((e^n_i)^*)_{i=1}^{4^n}$
 be the coordinate functionals satisfying $(e^{n}_i)^*(e^n_j) = \delta_{i,j}$ 
The mapping $f_n \colon \mathcal{L}_n \rightarrow X$ will be of the following form: \begin{equation} \label{eq: formoff_n}
f_n(a) = \sum_{i=1}^{4^n} (e^{n}_i)^*(f_n(a))e^n_i, \end{equation}
where  $(e^{n}_i)^*(f_n(a)) \in \{0,1\}$ and $\operatorname{supp}(f_n(a))= \{i \colon (e^{n}_i)^*(f_n(a)) = 1\}$ has size $|\operatorname{supp}(f_n(a))|= d(A,a)$. Note that $d(A,a)$  represents the  `depth' of $a$  in $\mathcal{L}_n$.

To define $f_n$ inductively, we suppose that $f_{n-1} \colon \mathcal{L}_{n-1} \rightarrow X$ has already  beeen defined to be of the form \eqref{eq: formoff_n} with $n$ replaced by $n-1$. 

Let $\rho \colon \mathcal{L}_{n-1} \rightarrow X$ be a `copy' of $f_{n-1}$ with $(e^{n-1}_i)_{i=1}^{4^{n-1}}$ replaced by $(e^n_i)_{i=1}^{4^{n-1}}$. The formal definition is as follows:
$$ \rho(a) = \sum_{i=1}^{4^{n-1}} (e^{n-1}_i)^*(f_{n-1}(a))e^n_i.$$ Similarly, let $\theta \colon \mathcal{L}_{n-1} \rightarrow X$ be a copy of  $f_{n-1}$ with $(e^{n-1}_i)_{i=1}^{4^{n-1}}$ replaced by $(e^n_i)_{i=4^{n-1}+1}^{2\cdot4^{n-1}}$. Formally,
$$ \theta(a) = \sum_{i=1}^{4^{n-1}} (e^{n-1}_i)^*(f_{n-1}(a))e^n_{4^{n-1}+i}.$$
Similarly,  let $\phi \colon \mathcal{L}_{n-1} \rightarrow X$ be a copy of  $f_{n-1}$ with $(e^{n-1}_i)_{i=1}^{4^{n-1}}$ replaced by $(e^n_i)_{i=2\cdot 4^{n-1}+1}^{3\cdot4^{n-1}}$. Formally,
$$ \phi(a) = \sum_{i=1}^{4^{n-1}} (e^{n-1}_i)^*(f_{n-1}(a))e^n_{2\cdot 4^{n-1}+i}.$$
Finally,  let $\sigma  \colon \mathcal{L}_{n-1} \rightarrow X$ be a copy of  $f_{n-1}$ with $(e^{n-1}_i)_{i=1}^{4^{n-1}}$ replaced by $(e^n_i)_{i=3\cdot 4^{n-1}+1}^{4^{n}}$. Formally,
$$ \sigma(a) = \sum_{i=1}^{4^{n-1}} (e^{n-1}_i)^*(f_{n-1}(a))e^n_{3\cdot 4^{n-1}+i}.$$
Recall that  $Y,C,D,E,F$ and $Z$ are `copies'  of $\mathcal{L}_{n-1}$.  Let $W$ be any one of these copies.  For    $a \in W$, let $\overline{a} \in D_{n-1} $   denote  the element of  $\mathcal{L}_{n-1}$ which corresponds to $a$. Now we define $f_n \colon 
\mathcal{L}_n  \rightarrow X$ as follows: \begin{equation*}  f_n(a) = \begin{cases} \rho(\overline{a}), & a \in Y\\
\sum_{i=1}^{4^{n-1}} e^n_i + \theta(\overline{a}), & a \in C\\
\sum_{i=1}^{4^{n-1}} e^n_i + \phi(\overline{a}), & a \in D\\
\sum_{i=1}^{2\cdot4^{n-1}} e^n_i + \phi(\overline{a}), & a \in E\\
\sum_{i=1}^{4^{n-1}} e^n_i +\sum_{i=2\cdot 4^{n-1}+1}^{3\cdot 4^{n-1}} e^n_i + \theta (\overline{a}), & a \in F\\
\sum_{i=1}^{3\cdot4^{n-1}} e^n_i + \sigma(\overline{a}), & a \in Z.
\end{cases}
 \end{equation*}
Note that    at the vertices $T,L,R$ and $B$, which connect the copies of $D_{n-1}$,   $f_{n}$ is defined twice, but both  definitions agree. So $f_n$ is well-defined.

Now we verify $\eqref{eq: bilipschitz}$. We begin with the right-hand inequality. If  $d(a,b)=1$, i.e.,  if $a$ and $b$ are adjacent vertices in $\mathcal{L}_n$, then it is clear from the definition that
$\|f_n(a) - f_n(b)\| \le 1$. Since $d$ is the shortest distance metric, the right-hand inequality follows at once from the triangle inequality in $X$.

We now turn to the left-hand inequality. If $a$ and $b$ belong to the same copy of $\mathcal{L}_{n-1}$ (either $Y,C,D,E,F$ or $Z$) then the left-hand inequality follows from the inductive hypothesis.  So suppose that they belong to different copies. There are several cases to consider.

\textbf{Case 1.}  Suppose that $a$ is `above' $b$ in $\mathcal{L}_n$.   Then $\operatorname{supp}(f_n(a) \subseteq \operatorname{supp}(f_n(b))$. Using Lemma~\ref{lem: J-convex},
\begin{align*}  \|f_n(b) - f_n(a)\|  &= \|\sum_{i \in \operatorname{supp}(f_n(b))\setminus \operatorname{supp}(f_n(a))} e^n_i\|\\
& \ge |\operatorname{supp}(f_n(b))| - |\operatorname{supp}(f_n(a))| - \varepsilon\\
& = d(a,b) - \varepsilon.
\end{align*}
\textbf{Case 2.}  Suppose $a \in C$, $b \in D$.  
\begin{align*} \|f_n(a) - f_n(b) \| &=\|\theta(\overline{a}) - \phi(\overline{b})\|\\
& \ge  |\operatorname{supp}(\theta(\overline{a}))| +  |\operatorname{supp}(\phi(\overline{b}))| - \varepsilon\\
\intertext{(by Lemma~\ref{lem: J-convex} since $\max \operatorname{supp}(\theta(a))< \min  \operatorname{supp}(\phi(b))$)}
& = d(T,a) + d(T,b) -\varepsilon\\
& = d(a,b) - \varepsilon.
\end{align*}
\textbf{Case 3.} Suppose $a \in C$, $b \in F$. Note that in this case $d(a,b) \le 2\cdot 4^{n-1}$. Hence
\begin{align*} f_n(a) - f_n(b) &= \theta(\overline{a}) - (\theta(\overline{b}) + \sum_{i=2 \cdot 4^{n-1}+1}^{3\cdot 4^{n-1}}e^n_i).
\end{align*} Note that $\theta(\overline{a}) - \theta(\overline{b}) =\sum_{i \in  A} \varepsilon_i e^n_i$, where $A \subseteq \{i \colon 4^{n-1}+1 \le i \le 2\cdot 4^{n-1}\}$ and $\varepsilon_i=\pm1$. Hence, by Lemma~\ref{lem: J-convex2},
\begin{align*}  \| f_n(a)-f_n(b)\| &=
\| \sum_{i \in  A} \varepsilon_i e^n_i  -  \sum_{i=2 \cdot 4^{n-1}+1}^{3\cdot 4^{n-1}}e^n_i\|\\
& \ge 4^{n-1} -   \varepsilon\\
&\ge \frac{1}{2} d(a,b) - \varepsilon.
\end{align*} \textbf{Case 4.} Suppose $a \in D, b \in E$.  This is similar to Case 3. Note that $d(a,b) \le 2 \cdot 4^{n-1}$.
\begin{align*}  f_n(a) - f_n(b) &= \phi(\overline{a}) - ( \sum_{i= 4^{n-1}+1}^{2\cdot 4^{n-1}}e^n_i + \phi(\overline{b}))\\
& =  (\phi(\overline{a}) - \phi(\overline{b})) - \sum_{i= 4^{n-1}+1}^{2\cdot 4^{n-1}}e^n_i. \end{align*}
 Note that $\phi(\overline{a}) - \phi(\overline{b}) =\sum_{i \in  A} \varepsilon_i e^n_i$, where $A \subseteq \{i \colon 2\cdot 4^{n-1}+1 \le i \le 3\cdot 4^{n-1}\}$ and $\varepsilon_i=\pm1$. Hence, by Lemma~\ref{lem: J-convex2},
\begin{align*}  \| f_n(a)-f_n(b)\| &=
\| \sum_{i \in  A} \varepsilon_i e^n_i  -  \sum_{i= 4^{n-1}+1}^{2\cdot 4^{n-1}}e^n_i\|\\
& \ge 4^{n-1} -   \varepsilon\\
&\ge \frac{1}{2} d(a,b) - \varepsilon.
\end{align*}  \textbf{Case 5.} Suppose  $a \in E, b \in F$. This is similar to Case 2. Note that
$$f_n(a) - f_n(b) = (\sum_{i= 4^{n-1}+1}^{2\cdot 4^{n-1}} e^n_i - \theta(\overline{b})) - (\sum_{i=2\cdot 4^{n-1}+1}^{3\cdot 4^{n-1}} e^n_i - \phi(\overline{a})).$$ Hence, by Lemma~\ref{lem: J-convex}, \begin{align*}
\|f_n(a) - f_n(b)\| &\ge (4^{n-1} - |\operatorname{supp}(\theta(\overline{b}))|) + (4^{n-1} - |\operatorname{supp}(\phi(\overline{a}))|) - \varepsilon\\
&= d(b,  B) + d(B, a) - \varepsilon\\
&= d(a,b) - \varepsilon. 
\end{align*}
\end{proof} 
\begin{remark}   The analogue of Theorem~\ref{thm: superreflexive} for $D_n$ is proved in \cite[Theorem 13.17, (13.26)]{P}  with the same distortion of $2+\varepsilon$.
 \end{remark}

We now prove a stronger result for $X = L_1[0,1]$.
\begin{theorem}  \label{thm: L1embedding} For each  $n \ge 1$, there exists a mapping $f_n \colon \mathcal{L}_n \rightarrow L_1[0,1]$ such that, for all $a,b \in \mathcal{L}_n$,  
\begin{equation}  \label{eq: bilipschitzL1} \frac{3}{4} d(a,b)  \le \|f_n(a)-f_n(b)\|_1 \le d(a,b).
\end{equation} \end{theorem} The proof requires the following elementary  lemma.
\begin{lemma}  \label{lem: elementary} For $0 \le s,t \le 1$,
$$1 + \min(s+t, 2-s-t) \le \frac{4}{3} (1+s+t-2st)$$
with equality if $s=t=1/2$. \end{lemma}
\begin{proof} First suppose $x:=s+t \le 1$. Then $\min(s+t, 2-s-t)=x$ and $st \le x^2/4$.  Hence  \begin{align*}
\frac{4}{3} (1+s+t-2st)-(1 + (s+t))&\ge \frac{4}{3} (1+x -\frac{x^2}{2})-1 -x\\
&= \frac{1}{3}  +\frac{x}{3} -  \frac{2x^2}{3}\\&\ge0.
\end{align*} The case $1 \le s+t$ is similar.
\end{proof} 
\begin{proof}[Proof of Theorem] Each $f_n$ will be of the following form: \begin{equation} \label{eq: defoff_n}
 f_n(a) = 4^n 1_{H_n(a)} \qquad(a \in \mathcal{L}_n), \end{equation}
where $H_n(a) \subseteq [0,1] $ and $|H_n(a)| = 4^{-n}d(A,a)$.  We begin with the base case $n=1$:
$$H_1(A) =\emptyset , H_1(T)= [0,1/4]; H_1(L) = [0,1/2]; $$$$H_1(R) = [0,1/4] \cup [1/2,3/4]; H_1(B)=[0,3/4]; H_1(U) = [0,1].$$ It is easily seen that $f_1$ is an isometry. 

For $n \ge 2$ the definition of $f_n$ is inductive. Suppose that $f_{n-1}$ has  been defined to be of the form  \eqref{eq: defoff_n}. Let  
 $\theta$ and $\phi$ be identically distributed copies of the mapping $a \mapsto H_{n-1}(a)$. Moreover, we require  $\theta$ and $\phi$ to be stochastically independent, i.e., 
$$|\theta(a)\cap \phi(b)| = |\theta(a)| ||\phi(b) |\qquad(a,b \in \mathcal{L}_{n-1}).$$
We use  $\theta$ and $\phi$ to define $H_n$ as follows:
\begin{equation*}  H_n(a) = \begin{cases}\frac{1}{4} \theta(\overline{a}), & a \in Y\\
[0,1/4] \cup(\frac{1}{4} + \frac{1}{4} \theta(\overline{a})), & a \in C\\
[0,1/4] \cup(\frac{1}{2} + \frac{1}{4} \phi(\overline{a})),  & a \in D\\
[0,1/2] \cup(\frac{1}{2} + \frac{1}{4} \theta(\overline{a})),  & a \in E\\
[0,1/4] \cup(\frac{1}{4} + \frac{1}{4} \phi(\overline{a})) \cup [\frac{1}{2},\frac{3}{4}] , & a \in F\\
[0,3/4] \cup(\frac{3}{4} + \frac{1}{4} \theta(\overline{a})), & a \in Z.
\end{cases}
 \end{equation*} The right-hand inequality of \eqref{eq: bilipschitzL1} follows as in the proof of Theorem~\ref{thm: superreflexive}. For the left-hand inequality we may assume that $a$ and $b$ belong to different copies of $\mathcal{L}_{n-1}$.

\textbf{Case 1.}  Suppose that $a$ is `above' $b$ in $\mathcal{L}_n$. Then $H_n(a) \subseteq H_n(b)$, so
$$ d(a,b) = 4^n(|H_n(b)| - |H_n(a)|) = \|f_n(a)-f_n(b)\|_1.$$

\textbf{Case 2.}  Suppose $a \in C$, $b \in D$.  Then \begin{align*}
\|f_n(a)-f_n(b)\|_1 
&= 4^{n-1}(|\theta(\overline{a})| + |\phi(\overline{b})|)\\
&= d(a,T) + d(b,T)\\
&= d(a,b).
\end{align*}

\textbf{Case 3.} Suppose $a \in C$, $b \in F$. Note that 
$$d(a,b) = 4^{n-1}(1+ \min(|\theta(\overline{a})| +|\phi(\overline{b})|, 2 - |\theta(\overline{a})| -|\phi(\overline{b})|)).$$
 Then \begin{align*}
\|f_n(a) - f_n(b)\|_1 &= 4^{n-1}(\|1_{\theta(\overline{a})}-1_{\phi(\overline{b})}\|_1+\|1_{0,1]}\|_1)\\
&=4^{n-1}(|\theta(\overline{a})| + |\phi(\overline{b})| -2|\theta(\overline{a})\cap \phi(\overline{b})| +1)\\
&=4^{n-1}(|\theta(\overline{a})| + |\phi(\overline{b})| -2|\theta(\overline{a})| |\phi(\overline{b})| +1)\\ \intertext{(since $\theta(\overline{a})$ and $\phi(\overline{b})$ are independent)}
&\ge 4^{n-1} \frac{3}{4}(1+ \min(|\theta(\overline{a})| +|\phi(\overline{b})|, 2 - |\theta(\overline{a})| -|\phi(\overline{b})|)) \\
\intertext{(from Lemma~\ref{lem: elementary}  with $s=|\theta(\overline{a})|$ and $t=|\phi(\overline{b})|$)}
&= \frac{3}{4}d(a,b).
\end{align*}

\textbf{Case 4.} Suppose $a \in D, b \in E$. This is essentially  the same as Case 3. As in Case 3, we obtain
$$\|f_n(a) - f_n(b)\|_1 \ge \frac{3}{4}d(a,b).$$

\textbf{Case 5.} Suppose  $a \in E, b \in F$. This is very similar to Case 2. Note that
 \begin{align*}
\|f_n(a)-f_n(b)\|_1 
&= 4^{n-1}((1-|\theta(\overline{a})| )+(1- |\phi(\overline{b})|))\\
&= d(a,B) + d(b,B)\\
&= d(a,b). \end{align*}
\end{proof} \begin{remark}  The analogue of Theorem~\ref{thm: L1embedding} for $D_n$ is proved in \cite[Theorem 5.1]{LR}  with the same distortion of $4/3$.
Moreover, it is remarked without proof that $4/3$ is the best constant for the distortion of  embeddings of $D_n$ as $n \rightarrow \infty$.  \cite[p. 359]{LR}.
In fact, we could not find  any  embedding of $D_2$ into $L_1[0,1]$ with distortion smaller than $4/3$ but we were not able to prove that $4/3$ is optimal. \end{remark}

 The next result shows that the distortion of any embedding of  $\mathcal{L}_2$  into $L_1[0,1]$ is at least $9/8$.
\begin{theorem} \label{thm: L_2dist}  Let $f \colon \mathcal{L}_2 \rightarrow L_1[0,1]$ satisfy 
$$d(a,b) \le \|f(a)-f(b)\|_1 \le c d(a,b).$$ 
Then $c \ge 9/8$.\end{theorem}

The proof uses the following result about hypermetric and negative type inequalities  from \cite{DL}.
\begin{theoremB}  \label{thm: negativetype}  \cite[Lemma 6.1.1]{DL} Let $(M\,\rho)$  be a finite metric space which is isometric to a subset of $L_1[0,1]$.   Then,  for all $k_i \in \mathbb{Z}$ ($1 \le i \le n$) such that $\sum_{i=1}^n k_i = 0$ (negative type inequalities) or $\sum_{i=1}^n k_i = 1$ (hypermetric inequalities), we have $$\sum_{1 \le i < j \le n} k_ik_j \rho(x_i, x_j) \le 0,$$
where $x_1,\dots,x_n$ are the distinct elements of $M$.
\end{theoremB} 
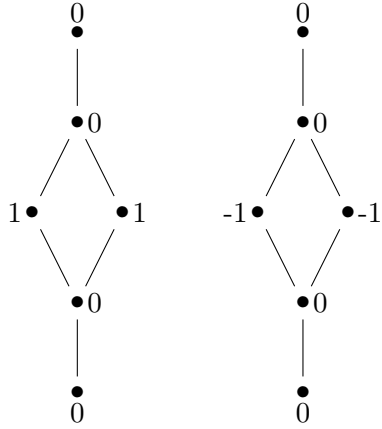
\begin{figure} \label{fig: L} \begin{center} {
\begin{tikzpicture}[scale=0.15]
 \node (A)  at (66,0) {$\bullet$};
\node[below] at (66, 0) {0};
\node (B) at (66,8) {$\bullet$}; \node[right] at (66,8) {0};
\node (C) at (62,16) {$\bullet$}; \node[left] at (62,16) {1};
\node (D) at (70,16) {$\bullet$}; \node[right] at (70,16) {1};
\node (E) at (66,24) {$\bullet$}; \node[right] at (66,24)  {0};
\node (F) at (66,32) {$\bullet$};
\node[above] at (66,32) {0};
\draw (A)  edge (B);
\draw (B)  edge (C); \draw  (B)  edge (D);\draw  (C)  edge (E); \draw  (D)  edge (E);
\draw  (E)  edge (F);
\node (A1) at  (86,0) {$\bullet$};
\node[below] at (86, 0) {0};
\node (B1) at (86,8) {$\bullet$}; \node[right] at (86,8) {0};
\node (C1) at (82,16) {$\bullet$}; \node[left] at (82,16) {-1};
\node (D1) at (90,16) {$\bullet$}; \node[right] at (90,16) {-1};
\node (E1) at (86,24) {$\bullet$}; \node[right] at (86,24)  {0};
\node (F1) at (86,32) {$\bullet$};
\node[above] at (86,32) {0};
\draw (A1)  edge (B1);
\draw (B1)  edge (C1); \draw  (B1)  edge (D1);\draw  (C1)  edge (E1); \draw  (D1)  edge (E1);
\draw  (E1)  edge (F1);
\end{tikzpicture}
}\end{center} \caption{Weights $P$ (left) and $N$ (right) for  $\mathcal{L}_1$} \label{figure: L1weighted} \end{figure}

 \begin{proof}[Proof of Theorem~\ref{thm: L_2dist}]  Consider the two choices of weights  for $\mathcal{L}_1$ indicated in Figure~\ref{figure: L1weighted} (each weight is  shown next to its corresponding vertex). 
Now define weights for $\mathcal{L}_2$ by assigning the $P$ weights to the $C $ and $F$ copies of $\mathcal{L}_1$, the $N$ weights to the $D$ and $E$ copies, and zero weights to the $Y$  and $Z$  copies. Let $(k_i)_{i=1}^{30}$ be the enumeration of these weights corresponding to some enumeration of the vertices of
 $\mathcal{L}_2$. Note that $\sum_{i=1}^{30}k_i=0$.
 By Theorem~B,
\begin{align*} 72 &= \sum_{i < j, k_ik_j>0} k_i k_j d(x_i,x_j)\\
&\le  \sum_{i < j, k_ik_j>0} k_i k_j \|f(x_i)-f(x_j)\|_1\\
&\le \sum_{i < j, k_ik_j<0} |k_i k_j| \|f(x_i)-f(x_j)\|_1\\
&\le c \sum_{i < j, k_ik_j<0} |k_i k_j |d(x_i,x_j)\\ &=64c.
\end{align*} So $c \ge 9/8$.
\end{proof} In a similar way  we can estimate the distortion of  metric embeddings of the diamond graph $D_2$ into $L_1[0,1]$.
 \begin{theorem} \label{thm: D_2dist} Let $f \colon D_2 \rightarrow L_1[0,1]$ satisfy 
$$d(a,b) \le \|f(a)-f(b)\|_1 \le c d(a,b).$$ 
Then $c \ge 5/4$.\end{theorem} \begin{proof} Consider the weights on $D_1$, denoted again $P$ and $N$, obtained from Figure~\ref{figure: L1weighted}
by removing the $A$ and $U$ vertices of 
$\mathcal{L}_1$.  Now define weights on $D_2$ by assigning $P$ to one pair of `opposite' copies of $D_1$ and $N$ to the other pair. Let $(k_i)_{i=1}^{12}$
be an enumeration of these weights corresponding to some enumeration of the vertices of
 $D_2$. Note that $\sum_{i=1}^{12}k_i=0$. Using Theorem~B as above yields
 $$40 =  \sum_{i < j, k_ik_j>0} k_i k_j d(x_i,x_j) \le  c \sum_{i < j, k_ik_j<0} |k_i k_j |d(x_i,x_j) =32c.$$
So  $c \ge 5/4$. \end{proof}
\begin{remark} 
A computer search revealed that $c = 5/4$ is the best estimate of the  lower bound for the distortion of $D_2$  that can be obtained from Theorem~B by considering all possible choices of $k_i$
in the range $-10 \le k_i \le 10$, and that $c = 9/8$ is the best that can be obtained for $\mathcal{L}_2$  by considering  all possible choices of $k_i$  in the range $-1 \le k_i \le 1$.
\end{remark} \begin{remark}  Since the proofs of  Theorems~\ref{thm: L_2dist} and \ref{thm: D_2dist} used only negative type inequalities, by \cite[Theorem 6.2.2]{DL} they  remain valid if $L_1[0,1]$ is replaced by $(\ell_2,\|\cdot\|_2^2)$. This is a stronger result as $L_1[0,1]$ is isometric to a subset of $(\ell_2,\|\cdot\|_2^2)$ (see e.g., \cite[p. 20]{O}).
\end{remark}

\end{document}